\documentclass[10pt]{article}
\usepackage[tbtags]{amsmath}
\usepackage{amssymb}
\usepackage{color}
\definecolor{c50}{rgb}{1,0,0}
\def\cL#1{\textcolor{c50}{#1}}
\def\cL#1{#1}
\allowdisplaybreaks[4] 
\newcommand{\nwc}{\newcommand}
\nwc{\COM}[1]{}

\nwc{\vs}[1]{\vskip #1 cm}

\newtheorem{theo}{Theorem}[section]
\newtheorem{sat}[theo]{Proposition}
\newtheorem{de}[theo]{Definition}
\newtheorem{lem}[theo]{Lemma}

\newtheorem{korr}[theo]{Corollary}
\newtheorem{remark}[theo]{Remark}
\newcommand{\nelem}[1]{{Lemma \ref{#1}}}
\newcommand{\neprop}[1]{{Proposition \ref{#1}}}
\newcommand{\netheo}[1]{{Theorem \ref{#1}}}

\newcommand{\kb}[1]{\boldsymbol{#1}}
\newcommand{\vk}[1]{\kb{#1}}

\def\FRE{\mbox{Fr\'{e}chet }}

\def\X{\vk{X}}

\newcommand{\ve}{\varepsilon}
\newcommand{\abs}[1]{\lvert #1 \rvert}
\newcommand{\Abs}[1]{ \Bigl \lvert #1 \Bigr \rvert}

\newcommand{\E}[1]{\mathbb{E}\{#1\}}
\newcommand{\pk}[1] {\mathbb{P} \left\{#1 \right\} }

\newcommand{\R}{\mathbb{R}}

\newcommand{\inr}{\in \R}

\newcommand{\ldot}{,\ldots,}

\newcommand{\limit}[1]{\lim_{#1 \to   \infty}}

\newcommand{\equaldis}{\stackrel{\mathcal{D}}{=}}

\newcommand{\BQN}{\begin{eqnarray}}
\newcommand{\EQN}{\end{eqnarray}}
\newcommand{\BQNY}{\begin{eqnarray*}}
\newcommand{\EQNY}{\end{eqnarray*}}
\newcommand{\BS}{\begin{sat}}
\newcommand{\ES}{\end{sat}}
\newcommand{\BRM}{\begin{remark}}
\newcommand{\ERM}{\end{remark}}

\newcommand{\BL}{\begin{lem}}
\newcommand{\EL}{\end{lem}}
\newcommand{\BT}{\begin{theo}}
\newcommand{\ET}{\end{theo}}
\newcommand{\BK}{\begin{korr}}
\newcommand{\EK}{\end{korr}}
\newcommand{\BD}{\begin{de}}
\newcommand{\ED}{\end{de}}

\newcommand{\QED}{\hfill $\Box$}

\newcommand{\IF}{\infty}

\def\fracl#1#2{\biggr( \frac{#1}{#2} \biggl) }

\newcommand{\prooftheo}[1]{ \textsc{Proof of Theorem} \ref{#1} }
\newcommand{\proofprop}[1]{\textsc{Proof of Proposition} \ref{#1}}
\newcommand{\prooflem}[1]{\textsc{Proof of Lemma} \ref{#1}}

\def\X{\vk{X}}

\def\lP{{\lambda(c,p)}}
\def\lP{\theta}
\def\lPY{\widetilde{\lambda(c,p)}}
\def\lPY{\widetilde{\theta}}
\def\Sp{S_p}

\begin{document}

\centerline{\Large Extremes of Aggregated Dirichlet Risks}

       \vskip 0.4 cm

\centerline{Enkelejd Hashorva\footnote{Faculty of Business and Economics (HEC Lausanne), University of Lausanne, 1015 Lausanne,  Switzerland, enkelejd.hashorva@unil.ch}}

\centerline{\today{}}

 {\bf Abstract:} The class of Dirichlet random vectors is central in numerous  probabilistic and statistical applications. 
The main result of this paper derives the exact tail asymptotics of the aggregated risk of powers of Dirichlet random vectors 
 when the radial component has df in the Gumbel or the Weibull max-domain of attraction. We present further results for the joint asymptotic independence and the max-sum equivalence.



\section{Introduction and Main Result}
Let $\X=(X_1 \ldot X_d)$ be a $d$-dimensional Dirichlet random vector with parameter $\vk{\alpha}=(\alpha_1 \ldot \alpha_d)\in (0,\IF)^d$ and radial component $R>0$ with some distribution function (df) $F$. By definition,  $\X$ has the stochastic representation
\BQN\label{stoc}
\X \equaldis \Bigl(R \frac{Y_1}{\sum_{i=1}^d Y_i}\ldot   R \frac{Y_d}{\sum_{i=1}^d Y_i}\Bigr)=:(R U_1 \ldot R U_d),
\EQN
where $\equaldis $ stands for equality of dfs and $Y_i,i\le d$ are independent random variables (rvs) such that $Y_i$ has
Gamma distribution  with parameters $\alpha_i$ and $1$ (in our notation the $Gamma(a,\lambda)$ distribution has probability density function (pdf) $\lambda^ax^{a-1}\exp(- \lambda x) /\Gamma(a)$ where $\Gamma(\cdot)$
is the Euler Gamma function).  Further $R, Y_1 \ldot Y_d$ and $\vk{U}=(U_1 \ldot U_d)$ are mutually independent.
 Basic distributional and asymptotic properties of Dirichlet random vectors are discussed in numerous contributions; see e.g., \cite{FanFan1990,KotBal2000,ChaSeg2009,McNNev2009,McNNev2010,Richter2011,Richter2012,MR2984355} and the references therein.

Clearly, for any $1 \le k< d$
$$ \sum_{i=1}^k X_i \equaldis R \sum_{i=1}^k U_i \equaldis R \frac{\sum_{i=1}^k Y_i}{\sum_{i=1}^d Y_i} \equaldis R B,$$
where $R$ and $B$ are independent, and $B$ has the Beta distribution with parameters $\sum_{i=1}^k\alpha_i$ and $ \sum_{k+1}^d \alpha_i$.
Hence the df of the total risk $\sum_{i=1}^k X_i$ can be directly calculated if $F$ is known. Clearly, when $k=2$ the above holds with $B$ almost surely equal to 1. Furthermore, if $F$ is in the Gumbel or the 
Weibull max-domain of attraction (MDA), then the tail asymptotics of $\sum_{i=1}^k X_i$ follows immediately by  Theorem 3.1 in \cite{Hasetal2010}. \\
In this paper we are concerned with the tail asymptotic behaviour of the aggregated risk $\Sp:=\sum_{i=1}^d \lambda_i X_i^p$ for some fixed 
constant $p>0$  and for given non-negative weights $\lambda_i,i\le d$. We shall assume first that  $\X$ with stochastic representation \eqref{stoc} has a radial component $R$ such that its df $F$ is in the Gumbel MDA, i.e., its survival function $\overline{F}=1-F$ satisfies for any $x\ge 0$
\BQN \label{eq:gumbel:w}
\overline{F}(u+x/w(u)) &\sim & \exp(-x) \overline{F}(u), \quad u\uparrow  x_F
 \EQN
for some positive  scaling function $w$ (here $x_F$ is the upper endpoint of $F$ and we abbreviate \eqref{eq:gumbel:w} as $F \in GMDA(w)$). We use in \eqref{eq:gumbel:w} the standard notation $\sim$ for the asymptotic equivalence of two real-valued functions.   For the sake of simplicity we shall assume hereafter that
$x_F=\IF$ or $x_F=1$. See \cite{Res1987,Faletal2010} for basic results concerned with the Gumbel MDA. Throughout in the following
\BQN\label{Lambda}
1=\lambda_1 = \cdots = \lambda_m \ge \lambda_{m+1} \ge \cdots \ge \lambda_d\ge 0
\EQN
are given weights  with $m \le d$ the multiplicity of $\lambda_1$. \cL{For} 
$p>1$ and $m< d$, it turns out that $\lambda_{m+1} \ldot \lambda_d$ do not influence the tail asymptotics of $\Sp= \sum_{i=1}^d \lambda_i X_i^p$,  which is however not the case if $p\in (0,1]$. Hereafter we set $\overline{\alpha}:= \sum_{i=1}^d \alpha_i$ with  $\alpha_i$'s being positive constants.\\
Our principal result below displays the exact asymptotics of the tail of $\Sp$, for any $p>0$.

\def\AAm{\overline{\alpha}_m}
\def\OAa{\overline{\alpha}}
\def\AM{\widehat \alpha}

\BT \label{Th1} Let $\X $ be a $d$-dimensional Dirichlet random vector with parameter $\vk{\alpha} $ and representation \eqref{stoc}. Suppose that \eqref{eq:gumbel:w} holds with $x_F\in \{1, \IF\}$ and some positive scaling function $w$.\\ 
a) If $p>1$, then 
\BQN\label{st:a}
\pk{\Sp> u^{p}}&\sim &\pk{\sum_{i=1}^m X_i^p > u^{p}}
\ \ \sim\ \ m^*\frac{\Gamma(\overline{\alpha})}{\Gamma(\AM) } ( u w(u))^{\AM -\overline{\alpha}} \overline{F}(u)
 , \quad u\uparrow  x_F, 
\EQN
where $\AM=\max_{1 \le i\le m} \alpha_i$, and $m^*$ is the number of elements of the index set $\{i\le m: \alpha_i= \AM\}$. \\
b) If $m< d$, then 
\BQN\label{st:b}
\pk{S_1 > u} &\sim&
\Biggl( \prod_{i=1}^{d-m} (1- \lambda_{m+i})^{-\alpha_{m+i}}\Biggr)  \frac{\Gamma(\overline{\alpha})}{ \Gamma(\sum_{i=1}^m \alpha_i)} (u w(u)) ^{- \sum_{i=1}^{d-m} \alpha_{m+i}} \overline{F}(u), \quad u\uparrow  x_F.
 \EQN
c) If $\lambda_i>0, i\le d$, then for any $p\in (0,1)$ we have
\BQN\label{st:c}
\pk{\Sp > \widetilde{\lambda_d} u^p} &\sim&
C_{\vk{\alpha}, d} (u w(u)) ^{ -(d-1)/2} \overline{F}( u), \quad u\uparrow  x_F,
 \EQN
with $C_{\vk{\alpha}, d}$ some positive constant and $\widetilde{\lambda_d}= \Bigl(\sum_{i=1}^d \lambda_i^{1/(1-p)}\Bigr)^{1-p}$.
\ET

{\bf Remarks}: a) An immediate consequence of \netheo{Th1} is that if $F$ is as therein, then the aggregated risk $S_p$ has df in the Gumbel MDA with scaling function $w_p(x)= x^{1/p-1} w(x^{1/p})/p$; see \cL{also} \neprop{Prop1} below. Consequently, in view of the properties of the scaling function $w$ (see e.g., \cL{p.143 in \cite{Embetal1997}}) we have assuming $x_F=\IF$ 
\BQNY \E{ S_p\lvert  S_p> VaR_{S_p}(b)} - VaR_{S_p}(b) &\sim& \frac{1}{w_p(VaR_{S_p}(b))}, \quad b \uparrow  1,
\EQNY
with $VaR_{S_p}(\tau)$ being the \cL{Value-at-Risk} of $S_p$ at $\tau \in (0,1)$, implying thus
\BQNY
\E{ S_p\lvert  S_p> VaR_{S_p}(b)} &\sim& VaR_{S_p}(b), \quad b \uparrow 1.
\EQNY
b)  For any df $F\in GMDA(w)$ with upper endpoint $x_F=\IF$,
the Davis-Resnick tail property is crucial, i.e., (see \cL{e.g., Proposition 1.1 in \cite{DavisR1988} and p. 113 in} \cite{HashExt12})
\BQN\label{DR}
\lim_{u\to \IF} (u w(u))^\mu \frac{ \overline{F}(c u)}{\overline{F}( u)}&=&0
\EQN
\cL{holds} for any $\mu\inr$ and $c>1$. Under the assumptions of statement $c)$ in \netheo{Th1} we have $\widetilde{\lambda_d}> \lambda_i, i\le d$.  It follows \cL{further} by \eqref{DR} that for $x_F=\IF, i\le d$ and $p\in (0,1]$
\BQN\label{st:c:2}
\limit{u} \frac{\pk{ \lambda_i X_i^p> u}}{\pk{\Sp > u}} &=& 0.
 \EQN
Consequently, each risk $\lambda_i X_i^p$ has a different asymptotic behaviour compared to $\Sp$. \\
$c)$ The convergence in \eqref{st:c:2} reveals a key  property of the Dirichlet dependence structure, namely the principle of a single big jump (see e.g., \cite{Foss2011} for details) applies if $p>1$. However, this principle does not apply when  $p\in (0,1]$, see 
\eqref{Prop22} below. An example which demonstrates this is furnished by taking $\X=(X_1 \ldot X_d)$ with independent components having unit exponential distribution, \cL{then} $\X$ is a Dirichlet random vector with its radial component having $Gamma(d,1)$ distribution.
Hence since also $\sum_{i=1}^d X_i$ has $Gamma(d,1)$ distribution if $p=1$, then 
$$ \lim_{u\to \IF} \frac{ \pk{\max_{1 \le i \le d} X_i^p> u}}{ \pk{ \Sp > u}}=0,$$
which is valid also for any $p\in (0,1)$.\\
$d)$ The tail asymptotic behaviour of $L_p$ type weighted norm $(\sum_{i=1}^d \lambda_i X_i^p)^{1/p}$ for various $\X$ has been considered by several authors; see e.g., \cite{PFatalov,MR2584909} and the references therein.

In the next section, we discuss our main result and present some  important extensions.  All the proof are relegated to Section 3 followed by an Appendix.

\section{Discussions and Extensions}
A canonical example of a $d$-dimensional Dirichlet random vector $\X$ is the so-called Kotz-Dirichlet random vector, with $X_i,i\le d$ independent such that $X_i$ has  $Gamma(\alpha_i,1)$ distribution with $\alpha_i>0,i\le d$; see e.g., \cite{MR2984355}.
Such a random vector has stochastic representation \eqref{stoc} with  $R$ having $Gamma(\overline{\alpha},1)$ distribution.
Hence for this particular example \netheo{Th1} gives the tail asymptotics of the sum of powers of independent Gamma rvs.\\ 
Note that for any $p>1$ the rv $X_i^p$
is a subexponential one (see e.g., \cite{Embetal1997} for \cL{the} definition and main properties), and therefore the statement a) in \netheo{Th1}
for this case can be directly checked to hold. When $p=1$, the claim of statement b) in \netheo{Th1} follows by Lemma 2.1 in \cite{Pakes2004}, whereas for $p\in (0,1)$ and $\alpha_i=\alpha>0, i\le d$ the claim in statement $c)$ of \netheo{Th1} is established by applying the result of \cite{Rootzen1}, which also gives the explicit formula for the constant $C_{\vk{\alpha},d}$ with $\vk{\alpha}=( \alpha \ldot \alpha)\in (0,\IF)^d$.

In the previous section we introduced the Dirichlet random vectors in the first quadrant. This restriction can be removed by introducing indicator rvs $I_1 \ldot I_d$ independent of $\X $ with stochastic representation \eqref{stoc}. If $\pk{I_i=1}=c_i=1- \pk{I_i=-1}, i\le d$, then
$$ \vk{Y}=(Y_1 \ldot Y_d) \equaldis (I_1 X_1^{1/p}  \ldot I_d X_d^{1/p}), \quad p>0$$
is referred to as a weighted $L_p$-Dirichlet random vector. 
 For simplicity, we assume here that $X_i, i\le d$ has the $Gamma(\alpha_i,1/p)$ distribution;
 the above extension allows us to include the Gaussian distribution in the class of $L_p$-Dirichlet random vectors. Indeed, if  $p=1/\alpha_1 =\cdots= 1/\alpha_d=2$ and  $I_1 \ldot I_d$ are mutually independent with mean 0, then $\vk{Y}$ is a $d$-dimensional Gaussian random vector if additionally $R^2$ is chi-square distributed with $d$ degrees of freedom. If  the df of $R$ is not specified in general, then $\vk{Y}$ is a spherical random vector (see the seminal contribution \cite{Cambanis} 
 for the main distributional properties).  We have thus $\sum_{i=1}^d \lambda_i \abs{Y_i}^2=\sum_{i=1}^d \lambda_i X_i$, and hence
for this particular case statement b) of \netheo{Th1} implies the claim of Theorem 3.1 in \cite{MR2584909}. 

In the sequel $\mathbb{B}_{a,b}$ stands for the Beta distribution with positive parameters $a$ and $b$, and $V\sim \mathbb{B}_{a,b}$ means that the rv  $V$ has the Beta distribution with parameters $a$ and $b$.\\ 
\COM{ To illustrate how this can be achieved, assume that $\lambda_{k+1}= \dots =\lambda_d=0$ for some $1< k<d$. We have
 $$ \sum_{i=1}^d \lambda_i X_i^p \equaldis R^p\sum_{i=1}^k \lambda_i U_i^p \equaldis (RB)^p \sum_{i=1}^k \lambda_i \tilde U_i^p=:(RB)^p Z_k ,$$
 where $B\sim \mathbb{B}_{\sum_{i=1}^k \alpha_i, \overline{\alpha}- \sum_{i=1}^k \alpha_i}$  and $(\tilde U_1 \ldot \tilde U_k)$
 is a standard Dirichlet random vector with parameter $\vk{\alpha}_k=(\alpha_1 \ldot \alpha_k)$. Furthermore, $R, B$ and $(\tilde U_1 \ldot \tilde U_k)$ are independent. Since further by the result in \eqref{eq:Lem1:2} (see \netheo{THconv} in Appendix) and statement $c)$ of \netheo{Th1} it follows that $B^p Z_m$ is also a regularly varying function at $\widetilde{\lambda_k}= (\sum_{i=1}^k \lambda_i^{1/(1-p)})^{1-p}$, then statement $c)$ in \netheo{Th1} can be easily adapted to include the case that $\lambda_i=0$ for some $i\le d$. \\
}

Concerning the Gumbel MDA assumption imposed on $F$ we first remark that under stronger assumptions on the scaling function $w$, 
namely $w$ is regularly varying at infinity,  then in view of \cite{HFarkasD13}, it follows that 
for any homogeneous function $h$ of order $p$, i.e., $h(t x_1 \ldot t x_d)= t^p h(x_1 \ldot x_d)$ holds for any $t>0$ and $(x_1 \ldot x_d) \cL{\inr^d}$ 
we have that $h(\X) \equaldis R^p h(\vk{U})$ has df in the Gumbel MDA. Using the  terminology of \cite{HashKorshPit} the rv $h(\X)$ can be referred to as the Dirichlet chaos. In the light of the findings of the aforementioned contribution, the exact asymptotics of the Dirichlet chaos can be derived.  In this paper we used a direct approach for the special case of aggregated risk. 

As mentioned in the Introduction the Davis-Resnick property of $F$ is crucial. In fact, if we assume that $\cL{\overline F}=1- F$ is rapidly varying at infinity, i.e., 
\eqref{DR} holds for $\mu=0$ and $c>1$, then for two Dirichlet random vectors $\X$ and $\vk{W}$ with corresponding radius $R$ and $R^*$ and parameter $\vk{\alpha}$,
 we obtain by  applying \nelem{XXH} in Appendix  
\BQN\label{uuu}
 \pk{ \sum_{i=1}^d \lambda_i X_i^p> u} &\sim & L(u) \pk{\sum_{i=1}^d \lambda_i W_i^p> u}, \quad u\to \IF ,
 \EQN
provided that $\overline F$ is rapidly varying at infinity and $\pk{R> u} \sim L(u) \pk{R^*> u}$ where $L(u)$ is some slowly varying function at infinity.

\subsection{Weibull MDA}
Instead of the Gumbel MDA assumption in \eqref{eq:gumbel:w} we shall suppose that $\overline{F}=1-F$ is regularly varying with index $\gamma\ge 0$ at the upper endpoint $x_F=1$, i.e., for any $t>0$
\BQN\label{eqWeib}
\frac{\overline{F}(1- tu)}{
\overline{F}(1- u)} &\sim& t^\gamma, \quad u \downarrow 0.
\EQN
For $\gamma>0$, the above assumption means that $F$ is in the MDA of the   Weibull distribution $\Psi_\gamma(x)=\exp(- \abs{x}^\gamma), x<0$. A canonical example of $F$ in the Weibull MDA is the case of the Beta distribution $\mathbb{B}_{a,b}$ where $\gamma=b$. Under \eqref{eqWeib} we can derive similar results to those in \netheo{Th1}. For simplicity we formulate only the claim of statement b) therein.

\BT\label{Th2}
Under the assumptions of statement b) in \netheo{Th1}, if further instead of \eqref{eq:gumbel:w} we suppose that 
the survival function $\overline{F}$ of $R$ satisfies \eqref{eqWeib} for some $\gamma\ge 0$, then 
\BQN
\pk{\Sp > 1- u}&\sim &
\Biggl(  \prod_{i=1}^{d-m} (1- \lambda_{m+i})^{-\alpha_{m+i}}\Biggr)
 \frac{\Gamma(\overline{\alpha})\Gamma(\gamma+1)}{\Gamma(\sum_{i=1}^{m}\alpha_i)\Gamma(\sum_{i=1}^{d-m}\alpha_{m+i}+ \gamma+1) }u^{-\sum_{i=1}^{d-m} \alpha_{m+i}} \overline{F}(1-u) \quad \quad 
 \EQN
holds as $u \downarrow 0.$
\ET
In the special case that $\alpha_1= \cdots = \alpha_d=1/2=1/p$ the claim of \netheo{Th2} agrees with that of Theorem 3.6 in \cite{MR2584909}.

\COM{
We discuss briefly the tail asymptotics of $\sum_{i=1}^d \Lambda_i X_i^p$ where $X_i$'s are as above but we take instead of $\lambda_i$'s some 
non-negative random weights $\Lambda_i,i\le d$ being independent of $\X$. Suppose that $\Lambda_i, i\le d$ has df with upper endpoint $\lambda_i$ being regularly varying at $\lambda_i$ with some index $\tau_i\ge 0$.
It is possible to extend \netheo{Th1} for such random weights, we state next the result for the case $p>1$ which can be shown along the same lines of the proof of \netheo{Th1}. The case $p\in (0,1]$ is more complicated and therefore omitted here.

\BT Under the assumptions of \netheo{Th1}, if $p>1$ and $1=\lambda_1= \cdots =\lambda_m> \lambda_{m+1} \ge \cdots \ge \lambda_d>0$ we have
\BQN
\pk{ \sum_{i=1}^d \Lambda_i X_i^p> u} &\sim &\sum_{i=1}^m \frac{ \Gamma( \overline{\alpha}- \alpha_i+ \tau_i +1)\Gamma(\overline{\alpha})} {\Gamma(\alpha_i)
\Gamma(\overline{\alpha}- \alpha_i+1)} (u w(u))^{- \overline{\alpha}+ \alpha_i} \pk{\Lambda_i> 1 - \frac{p}{u w(u)}} \overline{F}(u)
%
\EQN
as $u\uparrow x_F$.
\ET
}
\COM{
\BQN
\pk{ \sum_{i=1}^d \Lambda_i X_i^p> u^p} &\sim &\sum_{i=1}^n \pk{\Lambda_i X_i> u^p}\\
&\sim &\sum_{i=1}^n \pk{\Lambda_i^{1/p} R U_i> u}\\
&\sim &\sum_{i=1}^m \Gamma( \overline{\alpha}- \alpha_i+ \tau_i +1) \pk{\Lambda_i> 1 - \frac{p}{u w(u)}} \frac{\Gamma(\overline{\alpha})}{\Gamma(\alpha_i)
\Gamma(\overline{\alpha}- \alpha_i+1)} (u w(u))^{- \overline{\alpha}+ \alpha_i} \\
&\sim &\sum_{i=1}^m \frac{ \Gamma( \overline{\alpha}- \alpha_i+ \tau_i +1)\Gamma(\overline{\alpha})} {\Gamma(\alpha_i)
\Gamma(\overline{\alpha}- \alpha_i+1)} (u w(u))^{- \overline{\alpha}+ \alpha_i} \pk{\Lambda_i> 1 - \frac{p}{u w(u)}} \overline{F}(u)\\
\EQN
}

A specific of the Weibull MDA is that the upper endpoint $x_F$ of $F$ is necessarily finite. There is no possibility to convert $x_F$ to be infinite such that the transformed $\X$ is still 
a Dirichlet random vector. Therefore, the  result of this section cannot be retrieved by results available in the literature 
concerned with the aggregation of dependent unbounded risks dealt with for instance in  \cite{FossRich,MR3003975,Chavez}.

\subsection{Approximation by Max-Stable Distributions}
Next, we present an application of \netheo{Th1};  a similar application (omitted here) can be given using \netheo{Th2}.
Let $\vk{Y}=(Y_1 \ldot Y_d)$ be a random vector which is obtained by a linear transform of $(X_1^p \ldot X_d^p)$, i.e.,
for given constants $\lambda_{ij}, i,j\le d$
$$ \vk{Y} \equaldis (\sum_{i=1}^d\lambda_{i1}X_i^p \ldot \sum_{i=1}^d\lambda_{id}X_i^p).$$
We shall denote by $G$ the df of $\vk{Y}$, and $G_i$ is its $i$th marginal df. It is of interest to determine if $G$ is in the max-domain of attraction of some multivariate max-stable df $Q$, i.e., if there are constants $a_{ni}>0,b_{ni}\inr,i\le d,n\ge 1$ such that
\BQN\label{max-stable}
\lim_{n\to \IF}  \sup_{x_i \inr, 1 \le i\le d}
\Abs{  G^n( a_{n1} x_1 + b_{n1} \ldot a_{nd} x_d + b_{nd})- Q(x_1 \ldot x_d)}&=&0.
\EQN
Our next result shows that this is possible, if $F$ is in the Gumbel MDA.


\BS \label{Prop1}  Let $\lambda_{ij},i, j\le d$ be non-negative constants and denote by $A_j:=\{i\le d:\lambda_{ij}=1\}, j\le d$.
\cL{Suppose for $p\ge 1$ that} $A_j,j\le d$ is non-empty and $A_i\cap A_j$ has no elements for any pair $(i,j)$ of different indices, \cL{and for $p\in (0,1)$}  that $\sum_{i=1}^d \lambda_{ij}^{1/(1-p)}=1$ and $\lambda_{ij},i,j\le d$ are non-negative such that for any $i,j$ two different indices  $\lambda_{ik}\not=\lambda_{jk}$ for some $k\le d$.  Under the assumption of \netheo{Th1}, then for $a_{ni}=1/ w_p( b_{ni}),i\le d$ with $ b_{ni}=  G_i^{-1}( 1- 1/n), n\ge 1$ and $w_p(x)= x^{1/p-1} w(x^{1/p})/p,x>0$ we have that \eqref{max-stable} holds with
$Q(x_1 \ldot x_d)= \exp\Bigl(- \sum_{i=1}^d \exp(- x_i) \Bigr)$.
\COM{\BQN
\lim_{n\to \IF}  \sup_{x_i \inr, 1 \le i\le d}
\Abs{  \Bigl(G( a_{n1} x_1 + b_{n1} \ldot a_{nd} x_d + b_{nd})\Bigr)^n- \exp\Bigl(- \sum_{i=1}^d \exp(- x_i) \Bigr)}=0.
\EQN}
\ES

Clearly, the conditions in \neprop{Prop1} on $\lambda_{ij}$'s are satisfied if $\lambda_{ii}=1,i\le d$ and $\lambda_{ij}=0$ for all $i,j$ \cL{different} indices.
As in the proof of \neprop{Prop1} we have
\BQNY
\lim_{u\to \IF} \frac{ \pk{ X_i^p> u, X_j^p> u}}{\pk{ X_i^p> u}}&=&0.
\EQNY
Consequently, for the case $p>1$, by Bonferroni's inequality it follows that the sum and maximum of $\lambda_i X_i^p,i\le d$ are asymptotically equivalent, i.e.,  the principle of a single big jump holds. More precisely, if $x_F=\IF$ and $F \in GMDA(w)$, then for any $p> 1$
 \BQN \label{Prop22}
 \pk{\sum_{i=1}^d \lambda_i X_i^p> u} &\sim&  \pk{\max_{i\le d} \lambda_i X_i^p> u}, \quad u\to \IF.
\EQN

\subsection{Converse  Results}
So far we have assumed that the df of $R$ is in the Gumbel or Weibull MDA and then we showed that the same holds for the aggregated risk. Recall that we do not consider the 
case that $R$ has df in the \FRE MDA since the answer follows immediately by Breiman's lemma.\\
 At this point, the question on the 
validity of the converse results is natural. Namely,  
 if for some $\vk{\lambda}=(\lambda_1 \ldot \lambda_n)$ satisfying \eqref{Lambda}
the  aggregated risk $\Sp $ has df $G_{p,\vk{\lambda}}$ in the Gumbel MDA, then does also $F$ belong to the Gumbel MDA? 
Since in statistical applications, some observations might be missing, 
neither the radius $R$ 
\cL{nor} the total risk $\Sp$ can be observed, 
also of interest is  if $G_{p,\vk{\lambda}}$ belongs to the Gumbel MDA for some $\vk{\lambda}$ implies that $G_{p,\vk{\lambda}}$ is in the Gumbel MDA for any 
$\vk{\lambda}$ that satisfies \eqref{Lambda}. Note that when $F$ is in the Gumbel MDA, then $G_{p,\vk{\lambda}}$ is in the Gumbel MDA with scaling function $w_p(x)=x^{1/p-1}w(x^{1/p})/p$ for any $p\in (0,\IF)$.

We state next  the converse of \netheo{Th1} omitting the corresponding result for the Weibull MDA which can be derived by utilising the same idea.

\BT \label{BConv} Let $F,x_F,\X$ be as in \netheo{Th1}. If $\vk{\lambda}$ is a $d$-dimensional vector whose components satisfy  \eqref{Lambda}, 
 then $F \in GMDA(w)$ is equivalent with $G_{p,\vk{\lambda}} $ in the Gumbel MDA for some $\vk{\lambda}$ 
and some $p\in (0,\IF)$. Moreover, the latter assertion is equivalent with 
$G_{p,\vk{\lambda}}$ in the Gumbel MDA for any  $\vk{\lambda}$ and any  $p\in (0,\IF)$. 
\ET

Recent results concerning the asymptotics of products and converse results for the regularly varying case are derived in the deep contributions \cite{Jakobsen,JakobsenB}. Therefore, we omit the details for the case that $R$ has a regularly varying survival function at infinity. 


 
\section{Proofs}
We state first a lemma which is useful for the proof of \netheo{Th1}. In particular, the following lemma shows that in the bivariate setup \netheo{Th1} can be extended to include
some general bivariate random vectors which have similar dependence structure as the Dirichlet ones. In the sequel we say that $Z$ is regularly varying at $x_G $ with index $\tau\ge 0$  (we omit often 
the index $\tau$) if this is the case for its survival function $\overline G$.

\BL \label{A:00}  Let $B, X,Y$ be three non-negative rvs with upper endpoint\cL{s} $\omega_B=\omega_X=1, \omega_Y  \le 1$.\\
 a) If $\omega_Y < 1$ and $B^p X$ is regularly varying at 1 for some $p>1$, then 
for $\mathcal{S}_p:=B^p X+ (1-B)^p Y$ 
 \BQN\label{A:00:1}
 \pk{ \mathcal{S}_p> 1- u} \sim \pk{B^pX> 1- u}, \quad u\downarrow 0.
 \EQN
 b) Under the conditions of statement a), if \cL{further} $\omega_Y=1$ and $(1-B)^p Y$ is also regularly varying at 1, then
 \BQN\label{A:00:2}
 \pk{ \mathcal{S}_p > 1- u} \sim \pk{B^pX> 1- u}+ \pk{(1-B)^pY> 1- u}, \quad u\downarrow 0.
 \EQN
\COM{ c) For any $\lambda \in (-\IF, 1)$, if $B$ and $X$ are independent and regularly varying at 1 with  non-negative indexes $\alpha$ and $\gamma$, respectively,
 then
 \BQN\label{A:00:3}
 \pk{ B (X- \lambda)+ \lambda   > 1- u}  &\sim&   (1- \lambda)^{\alpha} \frac{\Gamma(\alpha+1)\Gamma(\gamma+1)}{\Gamma(\alpha+\gamma+1)}
 \pk{ B> 1- u} \pk{X> 1- u} , \quad u\downarrow 0.
 \EQN
}
 c) If $B$ has a continuous pdf $g$, then for any $c,\lambda$ positive  and $p\in (0,1)$
\BQNY
 \pk{ B^pc + \lambda (1-B)^p  > \lPY- u} & \sim  &2^{3/2}\frac{g(\lP)}{\sqrt{  h''(c,\lP)}}
\sqrt{u} , \quad  u\downarrow 0
\EQNY
holds with $h(c,\beta)= \beta^pc + \lambda (1-\beta)^p$ and $\lP=(\lambda/c)^{1/(p-1)}/(1+ (\lambda/c)^{1/(p-1)}),
\lPY=h(c,\lP)=(c^{1/(1-p)}+ \lambda^{1/(1-p)})^{p-1}.$\\
$d)$ Under the assumption and notation in statement $c)$ if further $X$ is regularly varying at $c:=\omega_X>0$  with index $\gamma> 0$,
then  for any $\lambda >0$ and $p\in (0,1)$ 
\BQNY
 \pk{ B^p X+ \lambda (1-B)^p  > \lPY- u}
& \sim  &\frac{\sqrt{2 \pi } g(\lP)}{\sqrt{  h''(c,\lP)}}
\frac{\Gamma(\gamma+ 1)}{\Gamma(\gamma+3/2)} \theta^{ - \gamma p} \sqrt{u} \pk{X > c- u},\quad u\downarrow 0,
\EQNY
provided that $B$ and $X$ are independent.
 \EL

\prooflem{A:00} a) For some $u>0$ sufficiently small, since $\omega_Y<1$, the event $\{ \mathcal{S}_p> 1- u\}$ is possible if
$B^p >1-u$ and $X> 1- u$ and thus in that case $(1-B)^pY=O(u^p)$. Hence
\BQNY \pk{ \mathcal{S}_p> 1- u} &\sim& \pk{B^pX> 1- u(1+o(1))}, \quad u\downarrow 0.
\EQNY
Thus the claim follows by the uniform convergence theorem for regularly varying function, see e.g.,  \cite{Embetal1997}.

b) As in the proof of a) the event  $\{ \mathcal{S}_p > 1- u\}$ is also possible if $B< u$ hence  $B^p X\le u^p$. Consequently
$$ \pk{ \mathcal{S}_p> 1- u} = \pk{B^pX> 1- u(1+o(1))}+ \pk{(1-B)^p Y> 1- u(1+o(1))} , \quad u\downarrow 0$$
and again the claim follows by the uniform convergence theorem for regularly varying function.

\COM{
c) In view of \eqref{againA} 
 \BQNY
 \pk{ B (X- \lambda)+ \lambda   > 1- u} 
 &\sim&   (1- \lambda)^{-\alpha} \frac{\Gamma(\alpha+1)\Gamma(\gamma+1)}{\Gamma(\alpha+\gamma+1)} \pk{ B> 1- u} \pk{X> 1- u} , \quad u\downarrow 0,
 \EQNY
since 
$$ \pk{ \frac{X - \lambda}{1- \lambda} > 1- u/(1-\lambda)}= 
\pk{X> 1- u}.$$
d) 
}
$c)$ First note that the unique maximum of the function $h(c, \beta)= \beta^pc + \lambda (1-\beta)^p$ for $\beta \in [0,1]$ is attained at
\BQN \label{lpt}
\lP=(\lambda/c)^{1/(p-1)}/(1+ (\lambda/c)^{1/(p-1)})
\EQN
 and we have thus $h'(c,\lP)=0$ and
$$\lPY= h(c,\lP)= \frac{\lambda}{(1+ (\lambda/c)^{1/(p-1)})^{p-1}}= c\lambda (c^{-1/(1-p)}+ \lambda^{-1/(1-p)})^{1-p} =(c^{1/(1-p)}+ \lambda^{1/(1-p)})^{1-p} .$$
Consequently, since $B$ has a continuous pdf $g$ we get that
for $\ve_u= \sqrt{2 u/h''(c,\lP)  }$
 \BQNY
 \pk{ B^p c+ \lambda (1-B)^p  > \lPY- u}  &\sim& \int_{\lP- \ve_u}^{\lP+ \ve_u} g(s)\, ds
  \sim  2^{3/2} \frac{g(\lP)}{\sqrt{h''(c,\lP)}}  \sqrt{u}
 \EQNY
as $u\downarrow 0$, hence the claim follows. \\
\cL{d}) Let $Q$ denote the df of $X$ and write $c>0$ for its upper endpoint.  Since $X$ is regularly varying at $c$ with index $\gamma>0$, then for any $t >0$
$$ \lim_{u\downarrow 0} \frac{\overline{Q}(c- tu)}{\overline{Q}(c- u)}= t^\gamma, \quad \overline{Q}=1- Q.$$
We proceed as above, but the choice of $\ve_u$ is different since we condition first on $X= c- tu$. 
Choosing  $\ve_u= \sqrt{\frac{2 u(1- {\lP}^p t)}{h''(c,\lP) } }$ with $\lP$ as in \eqref{lpt}, by the independence of $X$ and $B$ we may further write 
 \BQNY
\lefteqn{ \pk{ B^p X+ \lambda (1-B)^p  > \lPY- u} }\\
 &\sim&
 \int_{c-  u/\cL{\lP}^p}^c \int_{\lP- \ve_u}^{\lP+ \ve_u} g(s) \,ds d Q(t)\\ 
& \sim  & -2^{3/2} \frac{g(\lP)}{\sqrt{h''(c,\lP)}}  \overline{Q}(c- u)\sqrt{u}  \int_{0}^{1/ {\lP}^p}
\sqrt{ 1- {\lP}^p t}\, d \frac{Q(c- t u)}{\overline{Q}(c- u)}\\
& \sim  &2^{3/2} \frac{g(\lP)}{\sqrt{h''(c,\lP)}}  \overline{Q}(c- u)\sqrt{u} \gamma \int_{0}^{1/ {\lP}^p}
(1- {\lP}^p t)^{3/2- 1} t^{\gamma -1}\, dt\\
& \sim  &\sqrt{2 \pi } \frac{g(\lP)}{\sqrt{  h''(c,\lP)}}
\frac{\Gamma(\gamma+ 1)}{\Gamma(\gamma+3/2)}{\lP}^{- \gamma p}   \overline{Q}(c- u)\sqrt{u}
 \EQNY
 as $u\downarrow 0$, and thus the proof is complete. \QED

\prooftheo{Th1} In the sequel $\mathcal{B}_{\alpha,\beta}$ will denote a Beta rv with df $\mathbb{B}_{\alpha,\beta}$. 
Note that as $u\downarrow 0$
\BQN\label{bet} 
\pk{B^p> 1-u} \sim \frac{\Gamma(\alpha+\beta)}{\Gamma(\alpha)\Gamma(\beta)} \int_{1-u/p}^1 (1-x)^{\beta- 1} dx
\sim \frac{\Gamma(\alpha+\beta)}{p^\beta\Gamma(\alpha)\Gamma(\beta+1)} u^{\beta}.
\EQN
$a)$ Assume next that $m=1$, i.e., $1= \lambda_1> \lambda_2 \ge \cdots \ge \lambda_d\ge 0$. By the beta-independence splitting property of Dirichlet random vectors, we have the stochastic representation
\BQN\label{beta_ind}
 (U_1 \ldot U_d) \equaldis \Bigr( \mathcal{B}_{\alpha_1, \overline{\alpha}- \alpha_1},
(1- \mathcal{B}_{\alpha_1, \overline{\alpha}- \alpha_1})\tilde U_1 \ldot (1- \mathcal{B}_{\alpha_1, \overline{\alpha}- \alpha_1})\tilde U_{d-1}\Bigr),
\EQN
where $(\tilde U_1 \ldot \tilde U_{d-1})$ is a standard $(d-1)$-dimensional Dirichlet random vector with parameter $(\alpha_1 \ldot \alpha_{d-1})$ being independent of $\mathcal{B}_{\alpha_1, \overline{\alpha}- \alpha_1}$.
Consequently
$$ \sum_{i=1}^d \lambda_i U_i^p= \mathcal{B}_{\alpha_1, \overline{\alpha}- \alpha_1}^p+ \lambda (1-\mathcal{B}_{\alpha_1, \overline{\alpha}- \alpha_1})^p W,$$
where $\lambda \in (0,1)$ is some constant, $\mathcal{B}_{\alpha_1, \overline{\alpha}- \alpha_1}$ and $W$ are independent,
and $W$ has df with upper endpoint equal to 1. Applying statement a) of \nelem{A:00} we have as $u\downarrow 0$
$$
\pk{\sum_{i=1}^d \lambda_i U_i^p> 1- u} \sim \pk{U_1^p> 1- u}.
$$
\cL{Since further}
\BQN \label{repdp}
\sum_{i=1}^d \lambda_i X_i^p \equaldis R^p \sum_{i=1}^d \lambda_i U_i^p
\EQN
 and $R^p$ has df in the Gumbel MDA with scaling function $w_p(x)= x^{1/p-1} w(x^{1/p})/p$, see e.g., Lemma 5.2 in \cite{MR2678878}, 
the claim follows by applying \netheo{THconv}. Next, by repeating the above arguments, it follows that in the general case $1 \le m\le d$
$$
\pk{\sum_{i=1}^d \lambda_i U_i^p> 1- u} \sim \pk{\sum_{i=1}^m U_i^p> 1- u}, \quad u\downarrow 0.
$$
Since the case $m=1$ is shown above suppose that $m=2$. Again, by the beta-independence splitting property of Dirichlet random vectors
\BQN\label{C1} U_1^p+ U_2^p \equaldis \mathcal{B}_{\alpha_1, \overline{\alpha}- \alpha_1}^p
+ (1- \mathcal{B}_{\alpha_1, \overline{\alpha}- \alpha_1})^p \tilde U_2,
\EQN
with $\tilde U_2 \sim \mathbb{B}_{\alpha_2, \overline{\alpha}- \alpha_1- \alpha_2}$, provided that $d>2$. 
If
$d=2$, then we simply have
\BQN\label{C2} U_1^p+ U_2^p \equaldis \mathcal{B}_{\alpha_1, \alpha_2}^p+ (1- \mathcal{B}_{\alpha_1, \alpha_2})^p.
\EQN
In both cases, applying statement $b)$ and $c)$ of \nelem{A:00}
 we obtain
$$
\pk{U_1^p+U_2^p> 1- u} \sim  C_2 u^{\overline{\alpha}- \max(\alpha_1,\alpha_2)}, \quad u\downarrow 0,
$$
with $C_m \in (0,\IF),m\le d$. By induction on $m$ it follows that
$$
\pk{\sum_{i=1}^m U_i^p> 1- u} \sim  C_m u^{\overline{\alpha}- \max_{i\le m} \alpha_i}, \quad u\downarrow 0
$$
and further
$$
\pk{\sum_{i=1}^m U_i^p> 1- u} \sim  \pk{\sum_{1 \le i\le m: \alpha_i= \alpha^*_m} U_i^p> 1- u} , \quad u\downarrow 0,
$$
with $\alpha^*_m= \max_{i\le m }\alpha_i$. In order to simplify notation, assume that
$$\alpha_1=\alpha_i, \quad 2 \le i\le m^*\le m,$$
where $m^*$ denotes the number of elements in $\{1 \le i\le m: \alpha_i= \alpha^*_m\}$.
Suppose for simplicity that $m=m^*$ and consider next the case $m=2$. Clearly, if $d=2$, then by \eqref{C2} with $\alpha_1=\alpha_2$ and \nelem{A:00} it follows that (recall \eqref{bet})
\BQNY
\pk{\sum_{i=1}^{m^*} U_i^p> 1- u} &\sim&  m^*\pk{\mathcal{B}_{\alpha_1,\alpha_1}^p> 1- u}\\
 &\sim &   m^*\frac{\Gamma(2\alpha_1)}{\Gamma(\alpha_1)\Gamma(\alpha_1+1)}(u/p)^{\alpha_1} \\
 & \sim & m^*\pk{U_1^p> 1- u} , \quad u\downarrow 0.
 \EQNY 
For $d>2$ we consider the representation \eqref{C1}, where $\tilde U_2$ has Beta df with parameters $\alpha_2=\alpha_1, \overline{\alpha}- 2 \alpha_1>0$. In view of statement $b)$ and $c)$ of \nelem{A:00}
 we have
\BQNY
\pk{\sum_{i=1}^{m^*} U_i^p> 1- u} &\sim&  \pk{\mathcal{B}_{\alpha_1,\overline{\alpha}-\alpha_1}^p> 1- u}+
\pk{(1- \mathcal{B}_{\alpha_1,\overline{\alpha}-\alpha_1})^p \tilde U_2> 1- u}\\
 &\sim &   m^*\frac{\Gamma(\overline{\alpha})}{\Gamma(\alpha_1) \Gamma(\overline{\alpha}- \alpha_1+1)} (u/p)^{\overline{\alpha}-\alpha_1} , \quad u\downarrow 0.
 \EQNY
Using induction  and the above arguments, for any $m^*\ge 2$ we obtain
\BQNY
\pk{(\sum_{i=1}^{m^*} U_i^p)^{1/p}> 1- u} &\sim&  m^*\pk{U_1^p> 1- pu}\\
&\sim & m^*\pk{\mathcal{B}_{\alpha_1,\overline{\alpha}-\alpha_1}> 1- u} \\
&\sim&    m^*\frac{\Gamma(\overline{\alpha})}{\Gamma(\alpha_1) \Gamma(\overline{\alpha}- \alpha_1+1)} u^{\overline{\alpha}-\alpha_1} , \quad u\downarrow 0,
 \EQNY
hence the  claim follows by \netheo{THconv}.\\
$b)$ The case $m=d-1$ follows easily using the following representation
 $$\sum_{i=1}^m U_i+ \lambda_{m+1} U_{m+1} \equaldis
  B(1- \lambda_{m+1})+ \lambda_{m+1},$$
  where $B \equaldis  \mathcal{B}_{\sum_{i=1}^m \alpha_i,\alpha_{m+1}},$ and noting further that
\BQN\label{eqL0}
 \pk{ B(1- \lambda_{m+1})+ \lambda_{m+1}> 1- u}
 &\sim &(1- \lambda_{m+1})^{-\alpha_{m+1}} \frac{\Gamma(\overline{\alpha})}{\Gamma(\AAm)\Gamma(\overline{\alpha}-\AAm+1) }u^{\overline{\alpha}-\AAm}
 \EQN
as $u\downarrow 0$.  We consider next the case $m< d-1$. By the aggregation property of Dirichlet distributions and the beta-independence splitting property, we have
$$ \sum_{i=1}^{m+1} \lambda_i U_i \equaldis   B X+ \lambda_{m+1} (1- B),$$
where $B$ and $X$ are independent such that
$B \equaldis  \mathcal{B}_{  \overline{\alpha}- \alpha_{m+1},\alpha_{m+1} }$ and $  X\equaldis \mathcal{B}_{ \AAm,
\overline{\alpha}-\sum_{i=1}^{m+1} \alpha_i}. $ 
Consequently, \eqref{eq:Lem1:2} implies 
\BQNY
\pk{\sum_{i=1}^{m+1} \lambda_i U_i> 1- u} &\sim&
(1- \lambda_{m+1})^{-\alpha_{m+1}}   \frac{ \Gamma(\alpha_{m+1}+1) \Gamma(\overline{\alpha}-\sum_{i=1}^{m+1} \alpha_i+1)}
{\Gamma(\overline{\alpha}-\AAm+1)}\\
&& \times \pk{ \mathcal{B}_{  \overline{\alpha}- \alpha_{m+1},\alpha_{m+1} }> 1- u}
 \pk{\mathcal{B}_{ \AAm,
\overline{\alpha}-\sum_{i=1}^{m+1} \alpha_i} >  1- u} \\
&\sim&(1- \lambda_{m+1})^{-\alpha_{m+1}}
  \frac{ \Gamma(\overline{\alpha}) }{ \Gamma(\AAm)\Gamma(\overline{\alpha}-\AAm+1)} u^{\overline{\alpha}- \AAm}
 \EQNY
as $u\downarrow 0$.  Since \eqref{eqL0} holds also for $\lambda_{m+1}=0$,  repeating the above argument we have
\BQNY
\pk{\sum_{i=1}^{d} \lambda_i U_i> 1- u} &\sim& \Biggl(\prod_{i=1}^{d-m}(1- \lambda_{m+i})^{-\alpha_{m+i}}\Biggr)
\frac{\Gamma(\overline{\alpha})}{ \Gamma(\AAm) \Gamma(\overline{\alpha}- \AAm+1)}u^{\overline{\alpha}- \AAm}\\
&\sim& \Biggl(\prod_{i=1}^{d-m}(1- \lambda_{m+i})^{-\alpha_{m+i}}  \Biggr)\pk{ \sum_{i=1}^m U_i> 1- u}, \quad u\downarrow 0
 \EQNY
and hence the proof follows by applying again \netheo{THconv}.\\
$c)$ As above it suffices to determine  the tail asymptotics of $Z_d=\sum_{i=1}^d \lambda_i U_i^p$ at $\widetilde{\lambda_d}$
the upper endpoint of the df of $Z_d$. In view of \eqref{beta_ind} and statement $d)$ in \nelem{A:00} we have with 
$X:=\sum_{i=1}^{d-1}\lambda_i \tilde U_i$ being independent of $B \equaldis \mathcal{B}_{\overline{\alpha}- \alpha_d, \alpha_d}$
\BQNY
 \pk{Z_d> \widetilde{\lambda_d}- u}  &=&   \pk{B^p X + \lambda_d(1- B)^p > \widetilde{\lambda_d}- u}\\
& \cL{\sim} &\sqrt{2 \pi}  \frac{ g_{\overline{\alpha}- \alpha_d, \alpha_d}(\theta )} { \sqrt{ h''( \widetilde{\lambda_{d-1}},\theta )}} \frac{\Gamma(\alpha_d+1)}{\Gamma(\alpha_d+3/2)} \theta^{- \alpha_d p} \sqrt{u} \pk{X>
 \widetilde{\lambda_{d-1}}- u},  \quad u\downarrow 0,
 \EQNY
where $\widetilde{\lambda_{d-1}}$ is the upper endpoint of the df of $X$, $g_{\overline{\alpha}- \alpha_d, \alpha_d}$
is the pdf of  $B$ and
$$ \theta= \frac{\tau^{1/(p-1)}}{1+ \tau^{1/(p-1)}}, \quad \tau=\frac{\lambda_d}{\widetilde{\lambda_{d-1}}}.$$
From the proof of \nelem{A:00} we see that 
$$ \widetilde{\lambda_d}=  \frac{\lambda_d}{(1+ (\lambda_d/ \widetilde{\lambda_{d-1}})^{1/(p-1)})^{p-1}}
= ( \widetilde{ \lambda_{d-1}}^{1/(1-p)}+\lambda_d^{1/(1-p)})^{1-p},$$
hence
$$ \widetilde{\lambda_d}= \Bigl(\sum_{i=1}^d \lambda_i^{1/(1-p)}\Bigr)^{1- p}.$$
Note that above we used the fact that $X$ has a regularly varying survival function at $\widetilde{\lambda_{d-1}}$, which follows by induction. 
We remark further that $\widetilde{\lambda_d}$ is the attained maximum of the function $h(\beta_1 \ldot \beta_d)= \sum_{i=1}^d \lambda_i \beta_i^p$ for $\beta_i\in [0,1], i\le d$ 
satisfying  $\sum_{i=1}^d \beta_i=1$. Continuing, we obtain that
\BQNY
 \pk{Z_d> \widetilde{\lambda_d}- u}\sim   \pk{B^p X + \lambda_d(1- B)^p > \widetilde{\lambda_d}- u}\sim  \widetilde{C_d} u^{(d-1)/2},  \quad u\downarrow 0,
 \EQNY
with $\widetilde{C_d}$ a positive constant which can be calculated explicitly, and hence by \netheo{THconv}
\BQNY
 \pk{ \sum_{i=1}^d \lambda_i X_i^p> \widetilde{\lambda_d}u^p}
&=& \pk{ R (Z_d/\widetilde{\lambda_d})^{1/p} > u}\\
&\sim& \Gamma( (d-1)/2+1) \pk{  \frac{Z_d}{\widetilde{\lambda_d}} > 1- \frac{p}{u w(u)}}\pk{ R> u} \\
&\sim& \Gamma( (d+1)/2) \widetilde{C_d} \fracl{  p\widetilde{\lambda_d}}{u w(u)}^{(d-1)/2}\pk{ R> u} 
\EQNY
establishing the proof. \QED

\prooftheo{Th2} Applying \netheo{THconv} as in the proof of \netheo{Th1}, we obtain 
\BQNY
\pk{\sum_{i=1}^d \lambda_i X_i > 1- u}&\cL{=} &\pk{R \sum_{i=1}^d \lambda_i U_i  > 1- u } \notag\\ 
 &\sim & \frac{\Gamma(\sum_{i=1}^{d-m} \alpha_{m+i}+1)\Gamma(\gamma+1)}{\Gamma(\sum_{i=1}^{d-m} \alpha_{m+i}+ \gamma+1)} \pk{\sum_{i=1}^d \lambda_i U_i> 1- u} \pk{R> 1- u}  \notag\\
 &\sim&
 \prod_{i=1}^{d-m} (1- \lambda_{m+i})^{-\alpha_{m+i}}
 \frac{\Gamma(\sum_{i=1}^{d-m} \alpha_{m+i}+1)\Gamma(\gamma+1)}{\Gamma(\sum_{i=1}^{d-m}\alpha_{m+i}+ \gamma+1) }\\
&& \times  \frac{ \Gamma(\overline{\alpha})}{\Gamma(\sum_{i=1}^m \alpha_i) \Gamma(\sum_{i=1}^{d-m} \alpha_{m+i}+1)} u^{-\sum_{i=1}^{d-m} \alpha_{m+i}} \overline{F}(1-u)
\EQNY
as $u \downarrow 0$, hence the proof follows. \QED

\proofprop{Prop1}
In view of Theorem 3.1 and Lemma 5.2 in \cite{MR2678878}, it follows that $Y_j= \sum_{i=1}^d \lambda_{ij}X_i^p= R^p\sum_{i=1}^d \lambda_{ij}U_i^p$ has df in the Gumbel MDA with scaling function $w_p(x)= x^{1/p-1} w(x^{1/p})/p, x>0$, hence 
(see e.g., \cite{Faletal2010})
\BQNY
\lim_{n\to \IF}  \sup_{x_i \inr}
\Abs{  G^n_i( a_{ni} x_i + b_{n1})- \exp(- \exp(- x_i))}=0, \quad 1 \le i\le d.
\EQNY
Now by \cite{Res1987}, the claim follows if we show the pairwise asymptotic independence of $Y_i,Y_j$ for two different indices $i$ and $j$, i.e.,
\BQNY
\lim_{n\to \IF}  \frac{ \pk{ Y_i > b_{ni}, Y_j > b_{ni}}} {  \pk{Y_i> b_{ni}}} =0.
\EQNY
By the result of \netheo{Th1}, it follows that (see \cite{HashExt12})
$$ \limit{n} \frac{b_{ni}}{b_{n1}}=1, \quad 2 \le i \le d.$$
Clearly,
\BQNY
\frac{ \pk{ Y_i > b_{ni}, Y_j > b_{ni}}} {  \pk{Y_i> b_{ni}}} \le
\frac{ \pk{ Y_i+ Y_j  > 2b_{ni} (1+o(1))}} {  \pk{Y_i> b_{ni}}}
\EQNY
for all $n$ large. \cL{For} $p>1$, since by assumption $Y_i+Y_j= \sum_{k=1}^d 
 (\lambda_{ki}+ \lambda_{kj}) X_k^p$ with $\delta_k:=\lambda_{ki}+ \lambda_{kj}< 2$\cL{. A}pplying \netheo{Th1} we obtain
\BQNY
\frac{ \pk{ Y_i+ \lambda_j Y_j  > (1+\lambda_j) b_{ni} (1+o(1))}} {  \pk{Y_i> b_{ni}}}
&\to & 0, \quad n \to \IF,
\EQNY
which follows by the Davis-Resnick property mentioned in \eqref{DR}.  When $p=1$, the claim follows by statement $b)$ in \netheo{Th1} and \eqref{DR}. 
For $p\in (0,1)$, by the triangle inequality, and the assumption that  $\Bigl(\sum_{k=1}^d \lambda_{ki}^{\cL{q}}\Bigr)^{1/q}= \Bigl(\sum_{k=1}^d \lambda_{kj}^{q}\Bigr)^{1/q}=1$ with $q:=1/(1-p)$
, we have
$$ \widetilde{\delta_d}= \Bigl(\sum_{k=1}^d \delta_k^{q}\Bigr)^{1/q} < \Bigl(\sum_{k=1}^d \lambda_{ki}^{q}\Bigr)^{1/q}+ \Bigl(\sum_{k=1}^d \lambda_{kj}^{q}\Bigr)^{1/q}=2.$$
Hence statement c) of \netheo{Th1} and \eqref{DR} imply
 \BQNY
 \pk{ Y_i+ Y_j  > 2 b_{ni} (1+o(1))}&=& \pk{ Y_i+ Y_j  >  \widetilde{\delta_d} (2/\widetilde{\delta_d}) b_{ni} (1+o(1))}\\
&=&  o(\pk{Y_i> b_{ni}})
, \quad n\to \IF
 \EQNY
\COM{
hence by \eqref{uu} and \eqref{DR}
\BQNY
\frac{ \pk{ Y_i+ Y_j  > 2 b_{ni} (1+o(1))}} {  \pk{Y_i> b_{ni}}}&\to &0, \quad n \to \IF
\EQNY
}
and thus the claim  follows.
\QED

\COM{By the independence of $R$ and $U_i,i\le d$, we have using Theorem 3.1 in \cite{MR2678878} (recall $U_i$ has Beta distribution with parameters $\alpha_i$ and $\overline{\alpha}- \alpha_i$)
\BQN\label{uu}
 \pk{Y_i > u}= \pk{R U_i> u^{1/p}}&\sim& \Gamma(\overline{ \alpha}- \alpha_i+1) \pk{ U_i > 1- 1/(u^p w(u^p))} \overline{F}(u^{1/p})\notag\\
&\sim & \frac{ \Gamma(\overline{\alpha})}{\Gamma(\alpha_i) } \fracl{1}{ u^{1/p} w(u^{1/p})}^{\overline{\alpha}- \alpha_i}  \overline{F}(u^{1/p})
\EQN
as $u\to \IF$.
}

\prooftheo{BConv} In view of representation \eqref{repdp} and the tail behaviour of $\sum_{i=1}^d \lambda_i U_i^p$ found in the proof of \netheo{Th1},  the claim follows by applying 
\netheo{THconv} in Appendix. 
\QED

\section{Appendix} 
In \netheo{THconv} below we present results on the tail asymptotics of the products of two independent non-negative rvs. 
For its proof we need 
the next lemma, which is of some independent interest.

\BL \label{XXH} Let $S,S^*, Y,Y^*$ be \cL{four} independent positive rvs.  Let further 
$L$ be a slowly varying function at infinity and suppose that the dfs of $S$ and $S^*$ have upper endpoint equal to 1.\\
i) Assume that $\pk{S>x} \sim c \pk{S^*> x}$ \cL{as $x\uparrow 1$ for some $c\in (0,\IF)$.} If $Y$ has a rapidly varying survival function satisfying further $\pk{Y> u} \sim L(u) \pk{Y^*> u}$ as $u\to \IF$, then for any $\delta \in (0,1)$
\BQN
\pk{S Y > u} \sim c\pk{S^* Y> u} \sim \pk{SY>u, S > \delta } \sim  L(u)\pk{SY^*>u}, 
 \quad u\to \IF.
\EQN
ii) If $Y$ and $Y^*$ have dfs with upper endpoint equal to 1 and $\pk{Y> 1-1/u} \sim c^*\pk{Y^*> 1- 1/u}, c^*\in (0,\IF)$ as $u\to \IF$, then  we have
\BQN \label{26}
\pk{SY> 1- 1/u} \sim c^* \pk{S Y^*> 1- 1/u}, \quad u\to \IF.
\EQN
\EL

\prooflem{XXH} $i)$ Along the same lines of the proof of Lemma 1 in \cite{Chavez} for any $\delta \in (0,1)$  we have
\BQN\label{embT} \pk{SY> u} \sim \int_{\delta}^1 \pk{ Y> u/s} d G(s) = \pk{Y> u/\delta} \pk{S> \delta} + \int_{u}^{u/\delta} \pk{S> u/y}\, d F(y)
\EQN
as $u\to \IF$, where $F$ and $G$ are the dfs of $Y$ and $S$, respectively. Choosing $\delta$ close enough to 1 we obtain
$$ \pk{SY> u} \sim c\pk{Y> u/\delta} \pk{S^*> \delta} + c\int_{u}^{u/\delta} \pk{S^*> u/y}\, d F(y) \sim c\pk{S^*Y>u}$$
as $u\to \IF$. The other asymptotic equivalences are proved in \cite{HFarkasD13}, Lemma 4.1; the third
claim is due to Lemma A.3 in \cite{TanTsi2004}.\\
$ii)$ By the independence of $S, Y,Y^*$ for all $u$ and $G$ the df of $S$  we have
\BQNY
\pk{SY> 1- 1/u}&=&  \int_{1 - 1/u}^1 \pk{Y> (1- 1/u)/s} \, d G(s)\\
&\sim & c^*\int_{1 - 1/u}^1 \pk{Y^*> (1- 1/u)/s} \, d G(s)\\
\EQNY
as $u\to \IF$, hence the proof is complete.
\QED

{\bf Remark}: Let $S,S^*, Y$ be three non-negative independent rvs. Let $1$ be the upper endpoint of the dfs of $S$ and $S^*$ and suppose that the survival 
function of $Y$ is rapidly varying at infinity.  In view of Lemma 2 in \cite{Chavez} 
$$ \limit{u} \frac{\pk{SY> u}}{\pk{Y> u}}= \pk{S=1}, \quad \limit{u} \frac{\pk{S^*Y> u}}{\pk{Y> u}}= \pk{S^*=1},$$
hence \cL{if} $c=\pk{S=1}/\pk{S^*=1}> 0$, then  
\BQNY \pk{SY> u} &\sim& c \pk{S^* Y> u}, \quad u\to \IF.
\EQNY

\BT\label{THconv}
Let $S,Y$ be two independent non-negative rvs. Let $F$ and $H$ denote the dfs of $Y$ and $SY$, respectively. 
Suppose that for $L$ some slowly varying function at infinity  and some $\beta\ge 0$
\BQN\label{eq:betaTail:b}
\pk{S>1- 1/u}&\sim & L(u) u^{-\beta}, \quad u\to \IF.
\EQN
Assume further that $F$ has upper endpoint $x_F\in \{1, \IF\}$.\\
i) If $F \in GMDA(w)$, then  
\BQN\label{againA}
\pk{SY> u}&\sim & \Gamma(\beta+1) \pk{S> 1 - 1/(uw(u))} \pk{Y>u}, \quad u\uparrow x_F.
\EQN
Furthermore, if $\beta>0$ and $L(x)=L>0, \forall x>0$, then $H \in GMDA(w)$ if and only if 
 $F \in GMDA(w)$. \\
 ii) If  $F$ with $x_F=1$ satisfies \eqref{eqWeib} for some $\gamma \ge 0$, then for any $\lambda  \in (-\IF, 1)$
\BQN\label{eq:Lem1:2}
\pk{S(Y- \lambda)> 1-1/u} &\sim &  (1- \lambda)^\gamma \frac{\Gamma(\beta+1) \Gamma(\gamma+1)}{\Gamma(\beta+\gamma+1)}  
\pk{S> 1 - 1/u}\pk{Y> 1-1/u} , \quad u\to  \IF.
\EQN
Furthermore, if $\gamma >0, L(x)=L>0, \forall x>0$, then $F$ is in the Weibull  MDA of $\Psi_\gamma$ if and only if $H$ is in the Weibull MDA of $\Psi_{\beta+ \gamma}$. 
 \ET

\prooftheo{THconv} $i)$ Suppose that $x_F=\IF$. When  $S$ is beta distributed the claim follows from Theorem 4.1 in \cite{MR2678878}. 
Let us consider some general $S$ such that \eqref{eq:betaTail:b} holds. The claim in \eqref{againA} follows by Theorem 3.1 in  \cite{Hasetal2010}. Next, 
 we show that 
$H \in GMDA(w)$ implies $F \in GMDA(w)$. Since for any $\eta>1, u>0$ 
$$  \pk{S> 1/\eta}\pk{ Y> \eta u} = \pk{S> 1/\eta, Y> \eta u} \le  \pk{SY > u} \le \pk{Y>u}$$
and the fact that  $SY$  has df in the Gumbel MDA,  we conclude that both $SY$ and $Y$ have a rapidly varying survival function. 
If $L(t)=L>0, t>0$ and $\beta>0$, then for $\tilde S\equaldis \mathcal{B}_{a,\beta}$ with $a>0$ some arbitrary constant we find applying 
Theorem 3.1 in  \cite{Hasetal2010} that 
\BQNY \pk{\tilde SY> u} &\sim &\frac{\Gamma(a+\beta) }{L \Gamma(a)\Gamma(\beta+1)} \pk{SY> u}, \quad u\to \IF,
\EQNY
provided that $\tilde S$ is independent of $Y$. Hence  $\tilde S Y$ has df in the Gumbel MDA. It follows from Theorem 4.1 in \cite{MR2678878} that 
$Y$ has df in the Gumbel MDA with the same scaling function $w$ as $\tilde S Y$.  
\COM{We give below an alternative proof.  Condition \eqref{eq:betaTail:b} implies for any $\mu=\fracl{\Gamma(a+\beta) }{L\Gamma(a)\Gamma(\beta+1)}^{-1/\beta}$ with $a$ some positive constant 
\BQNY
\pk{S^{\mu}>1- u }\sim  \pk{S>1-    u (1+o(1)) /\mu } \sim L ( u/\mu)^\beta = \frac{\Gamma(a+\beta) }{\Gamma(a)\Gamma(\beta+1)} u^\beta, \quad u\downarrow 0.
\EQNY
Consequently if $\tilde S \eqauldis \mathcal{B}_{a,\beta}$, then we have (recall \eqref{bet})
$$\pk{S^{\mu}>1-u} \sim \pk{\tilde S>1-u}, \quad u\downarrow 0$$
and thus since $Y^{\mu}$ has also df in the Gumbel MDA (hence $Y^{\mu}$ has a rapidly varying survival function), then \nelem{XXH} implies taking  $\tilde S$ to be  independent of $Y$ 
$$\lim_{u\to \IF }  \frac{ \pk{(SY)^{\mu}>u}}{\pk{\tilde S Y^{\mu}>u}}=1.$$
This means that $\tilde S Y^{\mu}$ has df is in the Gumbel MDA. Hence by Theorem 4.1 in \cite{MR2678878} we have that $Y^{\mu}$ has df in the Gumbel MDA, 
which is equivalent with $Y$ has df in the Gumbel MDA. Hence the claim follows applying again Theorem 4.1 in \cite{MR2678878}. 
}
In view of \eqref{26} the case that  $x_F=1$ follows with similar arguments. 

$ii)$ The idea of the proof is the same as that of the proof of the statement $i)$ making further use of $ii)$ in \nelem{XXH}, Theorem 4.5 in \cite{MR2678878} and Theorem 3.1 in \cite{Hasetal2010}. \QED

\textbf{Acknowledgements}. I would like to thank the referees and an Editor for important comments and suggestions.  
Partial support from the Swiss National Science Foundation grants 200021-13478, 200021-140633/1  and the project RARE\cL{-}318984
 (an FP7  Marie Curie IRSES Fellowship) 	is kindly acknowledged.

\bibliographystyle{plain}
\bibliography{Arch1}

\def\polhk#1{\setbox0=\hbox{#1}{\ooalign{\hidewidth
  \lower1.5ex\hbox{`}\hidewidth\crcr\unhbox0}}}
\begin{thebibliography}{10}

\bibitem{Richter2012}
R.B. Arellano-Valle and W.-D. Richter.
\newblock On skewed continuous $l_{n,p}$-symmetric distributions.
\newblock {\em Chilean J. Stat.}, 3(2):193�--212, 2012.

\bibitem{MR2984355}
N.~Balakrishnan and E.~Hashorva.
\newblock Scale mixtures of {K}otz-{D}irichlet distributions.
\newblock {\em J. Multivariate Anal.}, 113:48--58, 2013.

\bibitem{Cambanis}
S.~Cambanis, S.~Huang, and G.~Simons.
\newblock On the theory of elliptically contoured distributions.
\newblock {\em J. Multivariate Anal.}, 11(3):368--385, 1981.

\bibitem{ChaSeg2009}
A.~Charpentier and J.~Segers.
\newblock Tails of multivariate {A}rchimedean copulas.
\newblock {\em J. Multivariate Anal.}, 100(7):1521--1537, 2009.

\bibitem{JakobsenB}
E.~Damek, T.~Mikosch, J.~Rosi{\'n}ski, and G.~Samorodnitsky.
\newblock General inverse problems for regular variation.
\newblock {\em Adv. Appl. Probab. to appear}, 2014.

\bibitem{DavisR1988}
R.~Davis and S.~Resnick.
\newblock Extremes of moving averages of random variables from the domain of
  attraction of the double exponential distribution.
\newblock {\em Stochastic Process. Appl.}, 30(1):41--68, 1988.

\bibitem{HFarkasD13}
K.~D{\polhk{e}}bicki, J.~Farkas, and E.~Hashorva.
\newblock Random scaling of {G}umbel risks.
\newblock {\em http://arxiv.org/abs/1312.7132}, 2013.

\bibitem{Chavez}
P.~Embrechts, E.~Hashorva, and T.~Mikosch.
\newblock Aggregation of log-linear risks.
\newblock {\em J. Appl. Probab.}, 2014, in press.

\bibitem{Embetal1997}
P.~Embrechts, C.~Kl{\"u}ppelberg, and T.~Mikosch.
\newblock {\em Modelling extremal events for insurance and finance}.
\newblock Springer-Verlag, Berlin, 1997.

\bibitem{Faletal2010}
M.~Falk, J.~H\"usler, and R.-D. Reiss.
\newblock Laws of {S}mall {N}umbers: {E}xtremes and {R}are {E}vents.
\newblock In {\em DMV Seminar}, volume~23. Birkh\"auser, Basel, third edition,
  2010.

\bibitem{FanFan1990}
K.T. Fang and B.Q. Fang.
\newblock Generalized symmetrized {D}irichlet distributions.
\newblock In {\em Statistical inference in elliptically contoured and related
  distributions}, pages 127--136. Allerton, New York, 1990.

\bibitem{Foss2011}
S.~Foss, D.~Korshunov, and S.~Zachary.
\newblock {\em An Introduction to Heavy-tailed and Subexponential
  Distributions}.
\newblock Springer-Verlag, New York, 2011.

\bibitem{FossRich}
S.~Foss and A.~Richards.
\newblock On sums of conditionally independent subexponential random variables.
\newblock {\em Math. Oper. Res.}, 35(1):102--119, 2010.

\bibitem{MR2584909}
E.~Hashorva.
\newblock Asymptotics of the norm of elliptical random vectors.
\newblock {\em J. Multivariate Anal.}, 101(4):926--935, 2010.

\bibitem{HashExt12}
E.~Hashorva.
\newblock Exact tail asymptotics in bivariate scale mixture models.
\newblock {\em Extremes}, 15(1):109--128, 2012.

\bibitem{MR3003975}
E.~Hashorva.
\newblock Exact tail asymptotics of aggregated parametrised risk.
\newblock {\em J. Math. Anal. Appl.}, 400(1):187--199, 2013.

\bibitem{HashKorshPit}
E.~Hashorva, D.~Korshunov, and V.I. Piterbarg.
\newblock Extremal behavior of {G}aussian chaos.
\newblock {\em arXiv:1307.5857v2}, 2013.

\bibitem{MR2678878}
E.~Hashorva and A.G. Pakes.
\newblock Tail asymptotics under beta random scaling.
\newblock {\em J. Math. Anal. Appl.}, 372(2):496--514, 2010.

\bibitem{Hasetal2010}
E.~Hashorva, A.G. Pakes, and Q.~Tang.
\newblock Asymptotics of random contractions.
\newblock {\em Insurance Math. Econom.}, 47(3):405--414, 2010.

\bibitem{Jakobsen}
M.~Jacobsen, T.~Mikosch, J.~Rosi{\'n}ski, and G.~Samorodnitsky.
\newblock Inverse problems for regular variation of linear filters, a
  cancellation property for {$\sigma$}-finite measures and identification of
  stable laws.
\newblock {\em Ann. Appl. Probab.}, 19(1):210--242, 2009.

\bibitem{KotBal2000}
S.~Kotz, N.~Balakrishnan, and N.~L. Johnson.
\newblock {\em Continuous multivariate distributions. {V}ol. 1}.
\newblock Wiley-Interscience, New York, second edition, 2000.

\bibitem{McNNev2009}
A.J. McNeil and J.~Ne{\v{s}}lehov{\'a}.
\newblock Multivariate {A}rchimedean copulas, {$d$}-monotone functions and
  {$l_1$}-norm symmetric distributions.
\newblock {\em Ann. Statist.}, 37(5B):3059--3097, 2009.

\bibitem{McNNev2010}
A.J. McNeil and J.~Ne{\v{s}}lehov{\'a}.
\newblock From {A}rchimedean to {L}iouville copulas.
\newblock {\em J. Multivariate Anal.}, 101(8):1772--1790, 2010.

\bibitem{Pakes2004}
A.G. Pakes.
\newblock Convolution equivalence and infinite divisibility.
\newblock {\em J. Appl. Prob.}, 41:407--424, 2004.

\bibitem{PFatalov}
V.I. Piterbarg and V.R. Fatalov.
\newblock The {L}aplace method for probability measures in {B}anach spaces.
\newblock {\em Uspekhi Mat. Nauk}, 50(6(306)):57--150, 1995.

\bibitem{Res1987}
S.I. Resnick.
\newblock {\em Extreme Values, Regular Variation, and Point Processes},
  volume~4 of {\em Applied Probability. A Series of the Applied Probability
  Trust}.
\newblock Springer-Verlag, New York, 1987.

\bibitem{Richter2011}
W.-D. Richter.
\newblock On skewed continuous $l_{n,p}$-symmetric distributions.
\newblock {\em Lithuanian Mathematical Journal}, 51:440--449, 2011.

\bibitem{Rootzen1}
H.~Rootz\'{e}n.
\newblock A ratio limit theorem for the tails of weighted sums.
\newblock {\em Ann. Proabb.}, 15:728--747, 1987.

\bibitem{TanTsi2004}
Q.~Tang and G.~Tsitsiashvili.
\newblock Finite- and infinite-time ruin probabilities in the presence of
  stochastic returns on investments.
\newblock {\em Adv. in Appl. Probab.}, 36(4):1278--1299, 2004.

\end{thebibliography}
\end{document}